\documentclass[journal]{IEEEtran}

\usepackage[noadjust]{cite}
\usepackage{graphicx,color,psfrag,subeqnarray}
\usepackage{amsmath}
\usepackage{amssymb}
\usepackage{amsthm,amssymb}
\usepackage{color}
\usepackage{listings}
\usepackage{multirow}
\usepackage{bm,array}
\usepackage{listings}
\usepackage[font=footnotesize]{caption}
\usepackage{subcaption}
\usepackage{bigstrut}

\usepackage{algorithm}
\usepackage[noend]{algpseudocode}

\makeatletter
\def\BState{\State\hskip-\ALG@thistlm}
\makeatother

\DeclareCaptionSubType[alph]{figure}
\newcommand\textlcsc[1]{\textsc{\MakeLowercase{#1}}}

\begin{document}

\title{Market-based Demand Response via Residential Plug-in Electric Vehicles in Smart Grids}

\author{Farshad~Rassaei,~\IEEEmembership{Student Member,~IEEE,} Wee-Seng~Soh,~\IEEEmembership{Member,~IEEE,}
	Kee-Chaing~Chua,~\IEEEmembership{Member,~IEEE}
	\thanks{}
	\thanks{The authors are with the Department of Electrical and Computer Engineering, National University of Singapore, Singapore 117583 (e-mail: farshad@u.nus.edu; weeseng@nus.edu.sg; eleckc@nus.edu.sg).}}

\maketitle

\begin{abstract}
Flexibility in power demand, diverse usage patterns and storage capability of plug-in electric vehicles (PEVs) grow the elasticity of residential electricity demand remarkably. This elasticity can be utilized to form the daily aggregated demand profile and/or alter instantaneous demand of a system wherein a large number of residential PEVs share one electricity retailer or an aggregator. In this paper, we propose a demand response (DR) technique to manage vehicle-to-grid (V2G) enabled PEVs' electricity assignments (charging and discharging) in order to reduce the overall electricity procurement costs for a retailer bidding to a two-settlement electricity market, i.e., a day-ahead (DA) and a spot or real-time (RT) market. 
We show that our approach is decentralized, scalable, fast converging and does not violate users' privacy. Extensive simulations show significant overall cost savings can be achieved for a retailer bidding to an operational electricity market by using the proposed algorithm. This technique becomes more needful when the power grid accommodates a large number of intermittent energy resources wherein RT demand altering is crucial due to more likely contingencies and hence more RT price fluctuations and even occurring the so-called \textit{black swan events}. Finally, such retailer could offer better deals to customers as well.
\end{abstract}

\begin{IEEEkeywords}
Black swan event, demand response (DR), electricity markets, plug-in electric vehicles (PEVs), residential load, retailer, smart grids, vehicle-to-grid (V2G).
\end{IEEEkeywords}

\IEEEpeerreviewmaketitle

\section{Introduction}
\IEEEPARstart{P}{lug-in Electric Vehicles} (PEVs) increase the elasticity of residential electricity demand profile substantially. In particular, when a retailer provides electricity to a large number of PEV owner customers, inherent flexibility in PEVs' demand, diversity in the usage patterns and energy storage capability of PEVs could be employed in order to save electricity costs. 

In addition, ambitious plans, aspiring incentives, subsidies and supports for introducing PEVs and plug-in hybrid electric vehicles (PHEVs) into the transport sector have been set in many countries \cite{travel2009sustainable}. The roadmap is advocating industries and governments to attain an overall  PEV/PHEV sales share of \textit{at least} 50\% for light duty vehicle (LDV) sales by 2050 worldwide \cite{222222}.

The future paradigm of electricity markets is a subject of question and debate as various changes are occurring that remould existing electricity generation, transmission and consumption formats. Massive integration of renewables, more efficient consumption thorough demand response (DR) techniques and transition of consumers to prosumers are some of such principal alterations \cite{strasser2015review}.



In a deregulated power grid, electricity retailers submit their demand bids to the wholesale market. For example, for a day-ahead (DA) market, these demand bids could often have both a desired power demand's quantity and a price component. This indicates that the retailer buys the specified power quantity, provided the market clearing price (MCP) is not higher than its offered price. This bidding process could be implemented in a few rounds to let the retailers modify and update their bids prior to the final clearance in the market. 


Pennsylvania–-Jersey-–Maryland (PJM) Interconnection runs the largest competitive wholesale electricity market on the globe \cite{ott2003experience}. The average prices of electricity per MWh for DA and RT markets which have been sold over year 2014 are very close, i.e., \$48.9539 and \$48.2063, respectively \cite{PJM}. However, this fact may misrepresent the nature of these two markets at the first glance. The details about hourly pricing data for DA and RT markets could be completely distinct and unpredictable at several days and/or hours. 

In fact, large spikes may be seen due to unexpected imbalances in supply and demand, for example, when a large production generator faces a black-out or temperatures are suddenly changing. Hence, the high uncertainty, particularly in the RT market, can remarkably impact the overall electricity procurement cost for a retailer. The spikes could occur more frequently when the whole grid is relying on numerous intermittent energy resources, e.g., wind farms and solar panels, where more RT price fluctuations and even so-called \textit{black swan events} \cite{taleb2010black} may be occurred.   

On the other hand, mean reversion theory tells us that prices and returns ultimately proceed back towards the mean or the average. This mean or average can be the historical average of the price or return or another sensible mean \cite{poterba1988mean}. In other words, it is not very likely that the unprecedented spikes keep on occurring and this lasts long.


\subsection{Summary of Technical Contributions}
In the paper, we adopt the flexibility and diversity in the power demand and availability time of PEVs from real world data. Then, we propose a cost minimizing algorithm suitable for the retailers dealing with existing operational markets using both offline demand shaping and online demand altering. The contributions of this article can be summarized as follows: 

\begin{itemize}
	\item We provide a fast converging and scalable DR technique which minimizes the overall electricity procurement cost for a retailer while preserving the privacy of individual users.     
	\item Our presented algorithm is capable to shape and alter the aggregated demand profile in response to DA and RT markets.     
	\item The approach offers a suitable mechanism for the retailer to decide when and how to respond to price fluctuations in the RT market.  
	\item We lay our algorithm in a game theory framework which has a Nash equilibrium (NE) guaranteeing that the proposed approach is yielding to all users' best turnovers.    
	\item Our presented results are based on an operational electricity market and available vehicle usage patterns.    
\end{itemize}

\subsection{Related Work}

Paper \cite{ott2003experience} describes the characteristics of the PJM DA and RT electricity markets. The author discusses that economic motivations make the DA and RT market prices converge in the bidding processes. Additionally, locational marginal pricing (LMP) based markets succour steady grid operations by using pertinent pricing signals to the retailers.

An overview of demand response (DR) and their various classifications in a deregulated electricity market is discussed in \cite{albadi2008summary}. The authors in \cite{6910332} compare different bidding rules in wholesale electricity markets when there exist PEVs and renewables' penetration in the power grid. 

The authors in \cite{vagropoulos2013optimal} present a two-stage stochastic optimization approach for an electric vehicle (EV) aggregator engaging in DA and regulation markets to reduce the energy cost by optimal bidding. Nevertheless, their proposed method impose some inconvenience on the customers and the aggregator should have access to private information of the EVs, e.g., arrival time, departure time and battery capacity. The same issue exists in the proposed method in \cite{vaya2015optimal}. In \cite{mohsenianoptimal}, the author discusses how a time-shiftable load, that may comprise of several time-shiftable subloads, can send demand bids to DA and RT markets to minimize its electricity procurement cost. Although this paper provides optimal closed-form solutions for bidding, they do not seem to be applicable for the residential sector wherein the retailer does not have detailed information about customers' preferences due to privacy concerns.  


In \cite{ISFAR}, we present a statistical modelling and a closed-form expression for PEVs' uncoordinated charging demand. Furthermore, in \cite{sgtrfar} we propose a decentralized demand shaping algorithm for a priori known desired demand profiles for the next day. Our paper \cite{farisgteu15} provides the idea of joint shaping and altering the demand briefly.

In \cite{kim_bidirectional_2013}, charging and discharging of PEVs are managed in order to maximize the social and individual welfare functions. However, in residential sector, defining appropriate utility and welfare functions for the individual users is very ambiguous.   


\begin{figure}[t]
	\centering
	\includegraphics[width=\linewidth,height=2.4in]{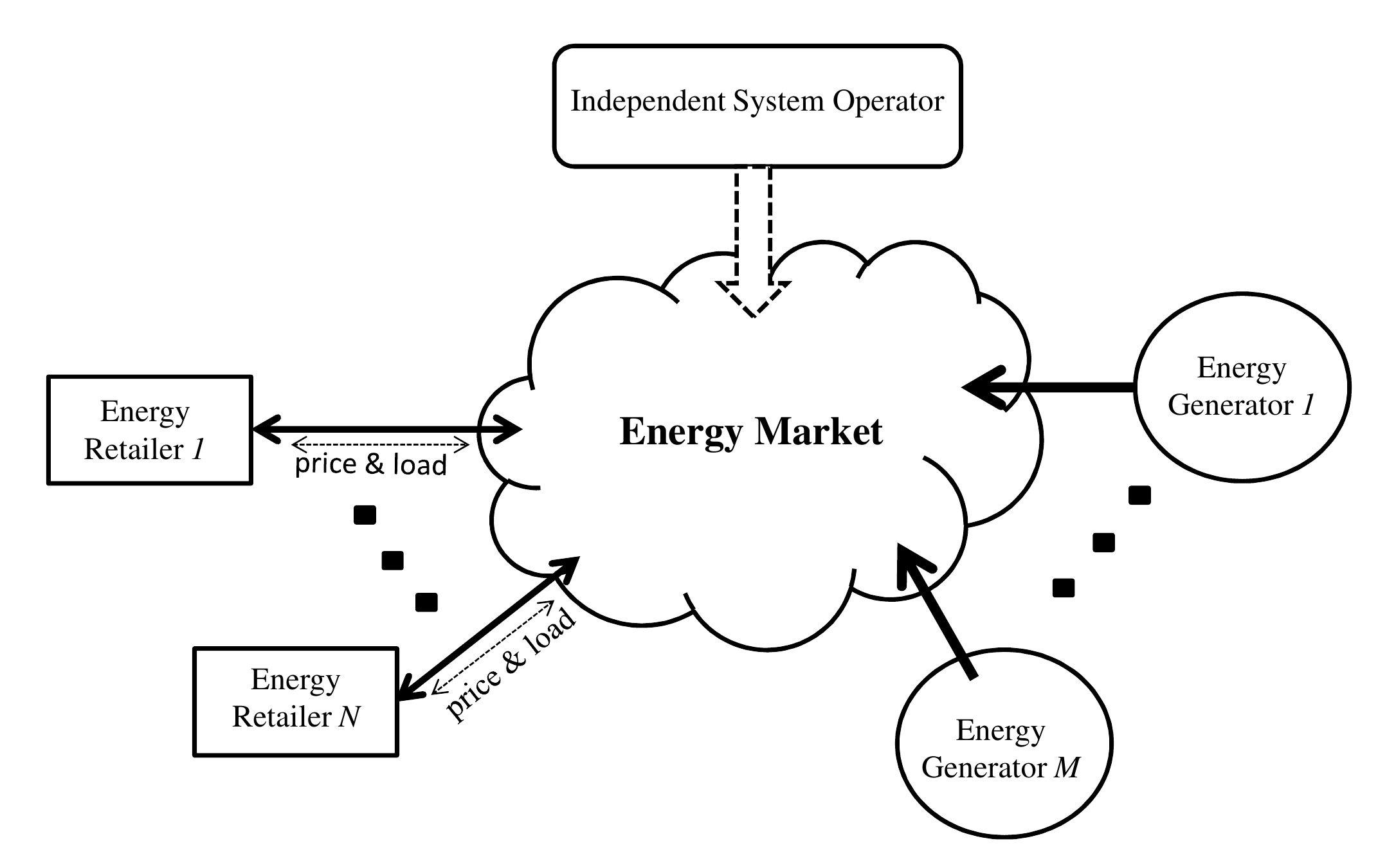} 
	\caption{A basic model of a smart grid comprised of energy market, multiple retailers and generators, and required communication infrastructure.} 
	\label{smodel}
\end{figure}
\begin{figure}[t]
	\centering
	\includegraphics[width=\linewidth,height=2.4in]{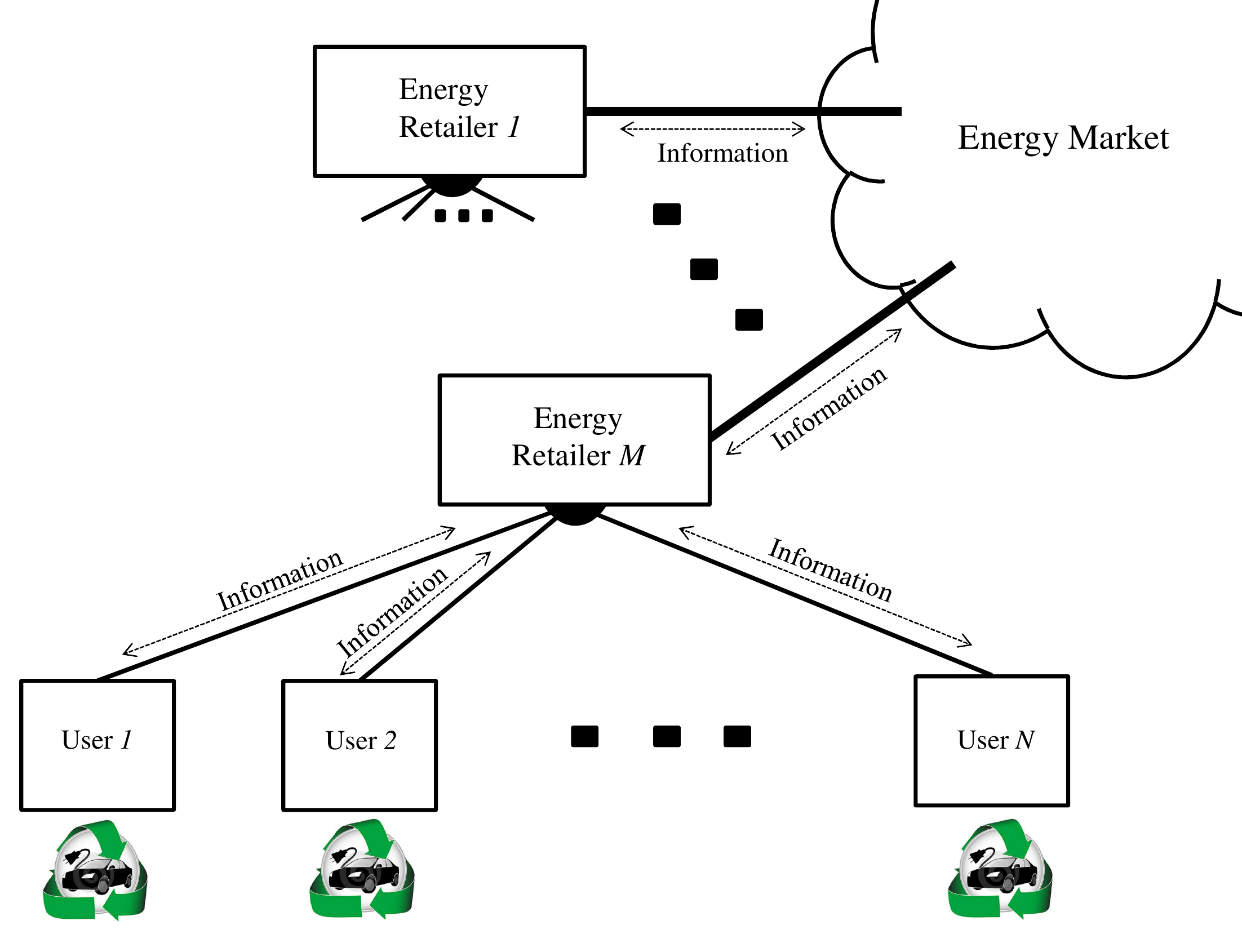} 	
	\caption{Multiple retailers serving their subscribed PEV owner residential users.} 
	\label{f1}
\end{figure}

%

In this paper, we provide a new architecture in which a joint demand shaping and altering algorithm is used to managing PEVs' electricity assignments (charging and discharging). The objective is minimizing the electricity procurement cost of a retailer bidding to two-settlement electricity markets, i.e., DA and RT markets. This algorithm includes both offline and online DR. We adopt PJM Interconnection pricing data \cite{PJM} for the year 2014 to evaluate our algorithm's results and efficacy.

\section{System Model} \label{SM}
In this section, we describe the underlying model of the power grid in this paper which contains the energy market, the electricity retailers, end users and necessary communication infrastructure. We articulate different parts of this model in the sequel.

Fig. \ref{smodel} shows a basic smart power grid model. Independent system operator (ISO) is supervising the market while generators and retailers deal electricity in the market. Here, we assume that retailers can also inject electricity back to the market and sell it. 

In a fully liberated electricity market, the bidding interactions between retailers and generators can be formulated and modelled as games, e.g., a \textit{Stackelberg game} \cite{von2010market}. In that case, the supervisory role of the ISO could be captured into the game's formulations and constraints. However, here, we consider an operational electricity market described by data.   

Fig. \ref{f1} represents multiple users share one electricity retailer or an aggregator. Different retailers exist in the power system and they compete to expand their market capacity by offering better deals to the customers. 

We assume that the users' overall load consists of two distinct types of load; typical household load which normally needs \textit{on-demand} power supply, e.g., air conditioning, lighting, cooking, refrigerator, etc, and PEV as a \textit{flexible} load. In this model, the dotted lines illustrates the underlying communication and information system while the solid lines show the power cabling infrastructure. 

We formed the pricing data for PJM market in 2013 into the following matrix for the DA market:

\begin{align} 
\textbf{P}^{DA}_{365,24} =
\begin{pmatrix}
p^{DA}_{1,1} & p^{DA}_{1,2} & \cdots & p^{DA}_{1,24} \\
p^{DA}_{2,1} & p^{DA}_{2,2} & \cdots & p^{DA}_{2,24} \\
\vdots  & \vdots  & \ddots & \vdots  \\
p^{DA}_{365,1} & p^{DA}_{365,2} & \cdots & p^{DA}_{365,24}
\end{pmatrix}.
\end{align} 

\noindent The same did we for the RT market. Then, the annual standard deviation for each hour of a day in PJM DA and RT markets is found as follows: 

\begin{align}
\nonumber \sigma(t)=\sqrt{\frac{1}{365}\sum_{i=1}^{365}(p^{DA/RT}_{i,t}-\frac{1}{365}\sum_{i=1}^{365}(p^{DA/RT}_{i,t})^2}, \\
\quad t=1, 2, \cdots, 24.
\label{stdf}
\end{align}

\begin{figure} 
	\centering
	\includegraphics[width=\columnwidth,height=2.4in]{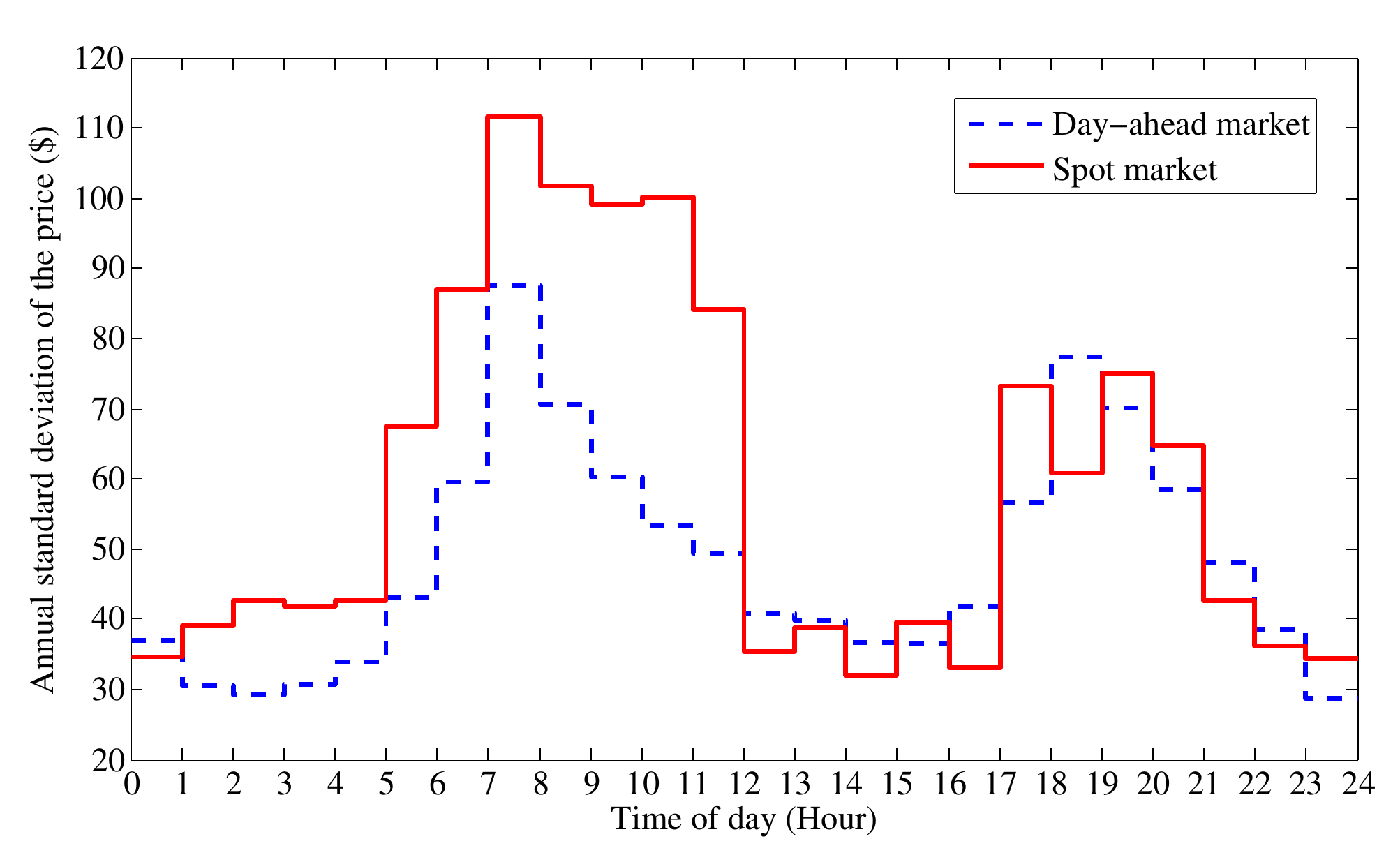}
	\caption[]{Annual standard deviation $\sigma(t)$ (\ref{stdf}) for each our of a day in PJM for DA and RT markets.}	 
	\label{std}
\end{figure}

Fig. \ref{std} shows the annual hourly standard deviation of the price for both DA and RT markets in 2013. As it can be observed and was expected, the spot market's prices can deviate much more for most hours of a day.   

We assume that a retailer prefers to shape its aggregated demand profile and emulate it to the profile purchased from DA market. Then, it is more able to minimize its demand from the RT market, which is much more prone to price volatility, to balance the load and accordingly lower the overall electricity procurement cost for each next day. This cost reduction makes the energy retailer capable to offer better deals to the customers in the form of cheaper pricing, random rewarding, promotions, etc, and expand its market capacity. 

In practice, the shaped aggregated profile for the next day does not exactly matches the retailer's purchased DA profile. Hence, the retailer often needs to reciprocate the load imbalances at each time slot of the following day by buying electricity from the RT market.

On the other hand, we should notice that residential DR is desired to be implemented such that users' privacy is not violated. Therefore, a DR practising retailer can expect more participation from the users if the actual DR is implemented in each user's house in a decentralized fashion according to our model. This also lessens the burden of heavy computations on a central unit and makes the algorithm more scalable.

\section{Analysis} \label{Sta}
In this section, we first formulate the electricity procurement cost for the retailer and then provide our proposed joint demand profile shaping and altering algorithm. 

In our analysis, we assume a daily energy assignment horizon and, without loss of generality, a time granularity of one hour. Let $(\textit{\textbf{l}}^d,\textit{\textbf{p}}^d)$ represents the pair of load $\textit{\textbf{l}}^d$ and price $\textit{\textbf{p}}^d$ vectors which has been cleared in the DA market. i.e.,

\begin{align}
\textit{\textbf{l}}^d\triangleq [l^d_{1},l^d_{2}\dotsc,l^d_{24}]^T,
\end{align} 
\begin{align}
\textit{\textbf{p}}^d\triangleq [p^d_{1},p^d_{2}\dotsc,p^d_{24}]^T,
\end{align}  

for which, the units of $l^d_i$ and $p^d_i$ are MWh and \$/MWh, respectively. Similarly, assume that $(\textit{\textbf{l}}^i,\textit{\textbf{p}}^r)$ represents the pair of load $\textit{\textbf{l}}^i$ and price $\textit{\textbf{p}}^r$ vectors which are the load imbalance and RT price vectors in the following day:

\begin{align}
\textit{\textbf{l}}^i\triangleq [l^i_{1},l^i_{2}\dotsc,l^i_{24}]^T,
\end{align} 
\begin{align}
\textit{\textbf{p}}^r\triangleq [p^r_{1},p^r_{2}\dotsc,p^r_{24}]^T.
\end{align}  

The values of the elements of these vectors will be only known to the retailer only at each time slot of the next day. 
Then, the overall electricity procurement cost for the next day can be formulated as follows: 
 
\begin{align}
C=\sum_{t=1}^{24}p^d_{t}l^d_{t}+\sum_{t=1}^{24}p^r_{t}l^i_{t},
\end{align}

here, $C$ is the overall energy procurement cost over the energy assignment horizon, i. e., 24 hours. 

First, given the purchased profile from the DA market by the retailer, i.e., $\textit{\textbf{l}}^d$, the users individually contribute to follow this demand profile by solving the sequential optimization problem \textbf{P1}. The objective is to minimize the correlation between each user's PEV energy assignment vector $\textit{\textbf{l}}_{\text{PEV},n}$ and its own inflexible demand vector $\textit{\textbf{l}}_{\text{A},n}$ plus the demand vector from other $N-1$ users $\textit{\textbf{l}}_{-n}$, and also to maximize the correlation between $\textit{\textbf{l}}_{\text{PEV},n}$ and the purchased DA load vector $\textit{\textbf{l}}^d$ (see Algorithm 1):

\begin{align}
\textbf{P1:} \quad \underset{\textit{\textbf{l}}_{\text{PEV},n}}{\text{minimize}}
\quad <\textit{\textbf{l}}_{\text{\text{PEV}},n},\textit{\textbf{l}}_{\text{A},n} + \textit{\textbf{l}}_{-n}-\textit{\textbf{l}}^d>,
\label{P1}
\end{align}\vspace{-1.2em}
\begin{align}
&\sum_{t=\alpha_{n}}^{\beta_{n}} l^{t}_{\text{PEV},n}=E_{\text{PEV},n},\\
&|l^{t}_{\text{PEV},n}| \leq p_{max},\quad \forall t\in\mathbb{T}^P_{\text{PEV},n}, \\ 
& l^{t}_{\text{PEV},n}=0,  \quad \forall t\notin\mathbb{T}^P_{\text{PEV},n},\\
& SOC^{t=\alpha}_{\text{PEV},n}+\sum_{k=\alpha+1}^{t}l'^k_{\text{PEV},n}\geq 0.2 \times C_{\text{PEV},n}, 
\forall t\in\mathbb{T}^P_{\text{PEV},n}.
\end{align}
\noindent in the above, $<\textit{\textbf{x}},\textit{\textbf{y}}>$ shows the inner product or correlation between vectors $\textit{\textbf{x}}$ and $\textit{\textbf{y}}$ and $\textit{\textbf{l}}_{\text{PEV},n}$ and $\textit{\textbf{l}}_{\text{A},n}$ show the energy assignment vectors for user n's PEV and the aggregated load from its household appliances, respectively. $E_{\text{PEV},n}$ is the $n^{\text{th}}$ user's required energy to be allocated to its PEV. Likewise, $\alpha_{n}$ and $\beta_{n}$ represent the arrival time and departure time of the PEV. Furthermore, $|l^{t}_{\text{PEV},n}| \leq p_{max}$ limits the maximum power that can be delivered to/from the PEV and $\mathbb{T}^P_{\text{PEV},n}$ represents the permissible charging time set or simply the set of time slots during the PEV's \textit{connection time} to the power grid. Additionally, $\textit{\textbf{l}}_{-n}$ is the aggregated power profile from other $N-1$ users described as follows: 

\begin{align}
\textit{\textbf{l}}_{-n}=\sum_{\substack{i \in {N} \\  i \neq n}}  (\textit{\textbf{l}}_{\text{PEV},i}+\textit{\textbf{l}}_{\text{A},i}). 
\label{others}
\end{align}
In (12), $C_{\text{PEV},n}$ is the total storage capacity of the user n's PEV and we assume that in case of employing V2G in the system, PEV's state of charge (SOC) should not fall below 20\% of that total capacity for emergency usage and in order to make sure that the adverse impacts on PEV's battery lifetime due to complete depletion are avoided.   

\begin{figure}[h]
	\centering
	\includegraphics[width=\columnwidth,height=2.2in]{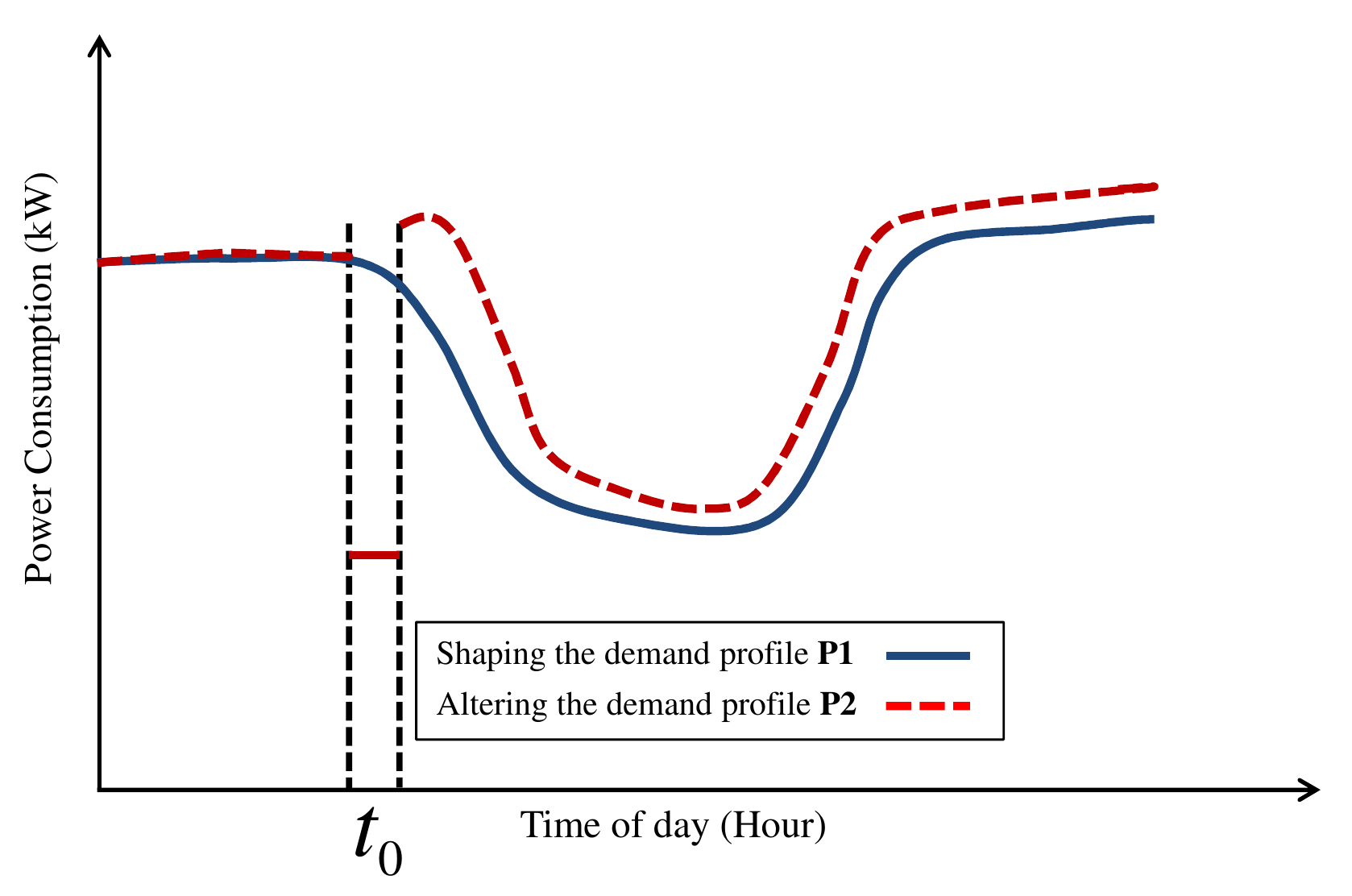}
	\caption{Altering the demand profile upon a contingency (lowering the demand at a time slot with unprecedented high RT price).} 
	\label{f4}
\end{figure}  
\begin{figure}[h] 
	\centering
	\includegraphics[width=\columnwidth,height=2.2in]{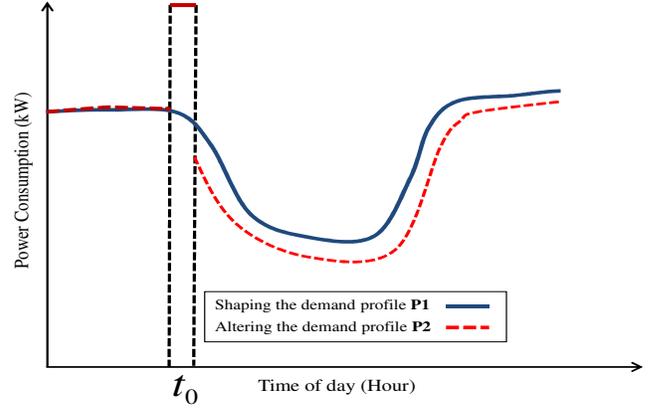}
	\caption{Altering the demand profile upon a contingency (lifting the demand at a time slot with unprecedented low RT price).} \vspace{0em}
	\label{f5}
\end{figure}
 
Second, knowing the fact that $(\textit{\textbf{l}}^i,\textit{\textbf{p}}^r)$ is unknown to the retailer a priori, at each time slot $t_0$ of the next day after getting this information, the retailer may decide to alter the previously shaped demand profile to minimize its RT electricity purchase to balance the load and even sell back some of its pre-purchased electricity from DA market if the RT price rises significantly due to the state of the RT market or contingencies.

As the chances for the price to remain that high during all the next few hours of the day is low, reshaping the load profile by lowering the electricity consumption at that time slot and purchasing electricity at the next time slots can yield to a lower cost in practice. This is also true for purchasing electricity at those time slots when price, unexpectedly, falls down significantly. The retailer may purchase extra electricity at those specific time slots (based on the overall storage capacity coming from connected PEVs) and reshape the demand for the next hours when the RT prices might be higher, c.f., Fig. \ref{f4} and Fig. \ref{f5}.

We should note that the retailer is assumed to be capable to indirectly employ the \textit{existing} flexibility (offered by each user's PEV). Nonetheless, the electricity consumption behaviours of the users (their PEVs' usage patterns) are not to be changed and hence the algorithm impose no restrictions on the users' comfort. Moreover, users' privacy is protected here.

In addition, we use a weighted moving average window as the threshold $\gamma_K(t)$ in Algorithm 1 based on which the retailer proceeds for RT demand altering:  
 
\begin{align}
\gamma_K(t)\triangleq \frac{1}{K}\sum\limits_{i=t-K}^{t} w_{i}p^r_{i}, 
\label{gama}
\end{align}

in which, $K$ is making the length of the window $K+1$ time slots and $w_{i}$ is defined as follows: 

\begin{align}
w_{i}\triangleq \frac{\sigma_(i)}{\sigma(t-K)}, \quad i=t-K, \cdots, t.
\end{align}

We should note that in the above formulation, for a negative $i$, it must be calculated in modulo 24. For instance, if $K=3$ and $t=1$, then $i=-2$ and $i \, mod\,24=22$. These weights can be calculated based on history data.

Then, the following linear multi objective programming (MOP) allows demand altering at time slot $t_0$ when $p_{t_0}^r$ is known and can be compared to $\gamma_K(t)$ along with pursuing the shape of the pre-purchased electricity from DA market (see Algorithm 1):

\begin{align}
\nonumber\textbf{P2:}  \quad \underset{\textit{\textbf{l}}'_{\text{PEV},n}}{\text{min}} \quad
b[\lambda <\textit{\textbf{l}}'_{\text{PEV,n}},\textit{\textbf{l}}'_{\text{A},n} + \textit{\textbf{l}}'_{-n}-\textit{\textbf{l}}^d>\\+ (1-\lambda) (l'^{t_0}_{-n}+l^{t_0}_{\text{A},n}+l'^{t_0}_{\text{PEV},n})],
\label{P2}
\end{align}\vspace{-1.2em}
\begin{align}
&[l^{'1}_{\text{PEV},n},\cdots,l^{'t_0-1}_{\text{PEV},n}]=[l^{1}_{\text{PEV},n},\cdots,l^{t_0-1}_{\text{PEV},n}]   ,\\
&\sum_{t=t_0}^{\beta_{n}} l'^{t}_{\text{PEV},n}=E_{\text{PEV},n}-\sum_{t=\alpha_{n}}^{t_0-1} l^{t}_{\text{PEV},n},  \\ 
&|l'^{t}_{\text{PEV},n}| \leq p_{max},\quad \forall t\in\mathbb{T}^P_{\text{PEV},n}, \\
& l'^{t}_{\text{PEV},n}=0, \quad \forall t\notin \mathbb{T}^P_{\text{PEV},n},\\ 
& SOC^{t=t_0-1}_{\text{PEV},n}+\sum_{k=t_0}^{t}l'^k_{\text{PEV},n}\geq 0.2 \times C_{\text{PEV},n}, 
\forall t\in\mathbb{T}^P_{\text{PEV},n}.
\end{align}
In the above, the value of b is either -1 or +1 according to Algorithm 1, $0\leq\lambda\leq1$ and $(1-\lambda)$ are the weights of the objective functions. In (21), $SOC^{t=t_0-1}_{\text{PEV},n}$ is the SOC of the user n's PEV before time slot $t_0$ and, similar to \textbf{P1}, we assume that in case of employing V2G in the system, the SOC should not fall below 20\% of the total PEV's storage capacity.   

We should note that when $\lambda=0$, ultimate possible flexibility that can be achieved from PEVs is obtained for a particular time slot. In this case, however, compliance with pre-purchased power profile from the DA market is sacrificed. On the other hand, when $\lambda=1$, there would be no altering in the profile. 

The value of $\lambda$ mainly depends on the price at time slot $t$, the storage capacity of connected PEVs at that time and the loss due to not following the priori shaped aggregated demand profile descriptively. In such cases, where it is not straightforward or impossible to formulate a function, using a fuzzy approach can help translate descriptive qualities into functions. An application of this technique is presented in \cite{6296552} for adjusting the step sizes in adaptive algorithms. 


The convergence criterion in Algorithm 1 can be simply assumed as a desired number of iterations for updating all users' demand profiles or it can be set and subjected to some predetermined error function, e.g., a desired mean square error (MSE) between two subsequent iterations of achieving aggregated demand profiles. Furthermore, similarly, as discussed in \cite{sgtrfar}, a convergence is guaranteed to be obtained and users' contribution can be modelled as a \textit{cooperative game} with \textit{complete information} wherein a Nash equilibrium exists \cite{sgtrfar}.      


\begin{algorithm}[t]
	\begin{algorithmic}[1]
		\caption{Offline \& Online Demand Response}\label{demshape}
		\State Each user initializes its respective load profile over the assignment horizon based on its power demands, i.e., $\textit{\textbf{l}}_{n}$ for $n=1,\dotsc, N$.
		\State All $N$ users send their initialized load profiles to the retailer.
		\While {\textit{not reaching convergence}}
		\For {$n=1$ to $N$}
		\State  The retailer calculates the state information $\textit{\textbf{l}}_{-n}$ according to (\ref{others}) for user $n$.
		\State The retailer sends $(\textit{\textbf{l}}_{-n}-\textit{\textbf{l}}^d)$ to user $n$.
		\State User $n$ solves problem \textbf{P1} and updates its load profile $\textit{\textbf{l}}_{n}$.
		\State User $n$ sends back the new demand profile to the retailer. 
		\State The retailer updates $\textit{\textbf{l}}_{n}$.
		\EndFor 
		\State \textbf{end for}
		\EndWhile 
		\State \textbf{end while}
		\For {$t=1$ to $24$}
		\State The retailer receives information from RT market, i.e., $\textit{p}^r_t$.
		\If {$(\textit{p}^r_{t} \neq \gamma_K(t) )$} 
		\If {$(\textit{p}^r_{t} > \gamma_K(t) )$}    
		\State $b=+1$.
		\Else 
		\State $b=-1$.
		\EndIf
		\State \textbf{end if}
		\State The retailer proceeds for demand altering.
		\While {\textit{not reaching convergence}}
		\For {$n=1$ to $N$}
		\State The retailer sends demand altering signal at time slot $t$ to user $n$.
		\State User $n$ solves problem \textbf{P2} and updates its load profile $\textit{\textbf{l}}'_{n}$.
		\State User $n$ sends back the new demand profile to the retailer. 
		\State The retailer updates $\textit{\textbf{l}}'_{n}$.
		\State  The retailer calculates the state information $\textit{\textbf{l}}'_{-n}$ according to (\ref{others}) for user $n$.
		\EndFor 
		\State \textbf{end for}
		\EndWhile 
		\State \textbf{end while}
		\EndIf 
		\State \textbf{end if}
		\EndFor 
		\State \textbf{end for}
	\end{algorithmic}
\end{algorithm}

\section{Simulation Results} \label{SR}
In this section, we illustrate the efficacy of the proposed algorithm in the previous section through computer simulations. In our simulations, the number of users (residential electricity consumers with PEVs), $N$, is 1,000 and the optimization horizon is considered to be a day, i.e., 24 hours for a DA programming scenario with a time granularity of one hour, i.e., duration of each time slot is one hour. Also, $K$ is considered to be 3 in (\ref{gama}). 

For the PEVs usage patterns, i.e., the arrival times, departure times and vehicles' energy demands, our data and distributions are based on 2009 NHTS data \cite{NHTS}. We considered new standard outlets, NEMA 5-15, with 1.8 kW power output. Furthermore, SOC for each PEV at the arrival time is as follows in percentage points: 

\begin{equation}
SOC^{t=\alpha}_{\text{PEV},n}=100 \times (1-\frac{E_{\text{PEV},n}}{24}).
\end{equation} 

In other words, we assumed that PEVs are fully charged by their respective next departure time. Additionally, we considered 24 kWh energy storage capacity for PEVs according to Nissan Leaf model \cite{leaf}. We adopted the PJM interconnection electricity market pricing data for both DA and RT markets in year 2014 \cite{PJM}. 

\begin{figure} 
	\centering
		\includegraphics[width=\columnwidth,height=2.4in]{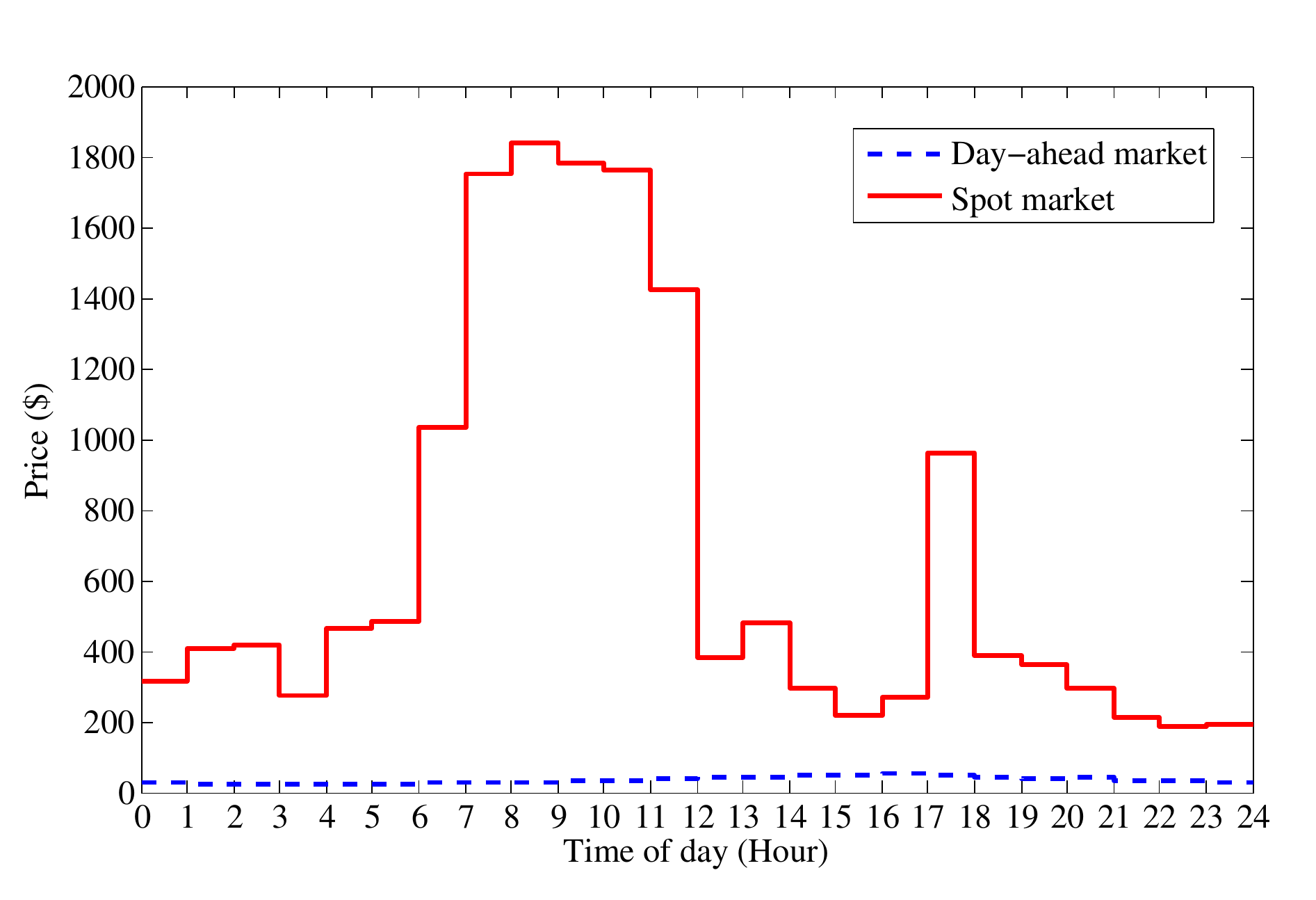}
	\caption[]{DA and RT prices for 9 March 2014 in PJM Interconnection electricity market (the number of days with such black swan behaviour in 2014 is quite considerable).} 
	\label{price}
\end{figure}

\begin{figure}[t]
	\centering
	\includegraphics[width=\columnwidth,height=2.4in]{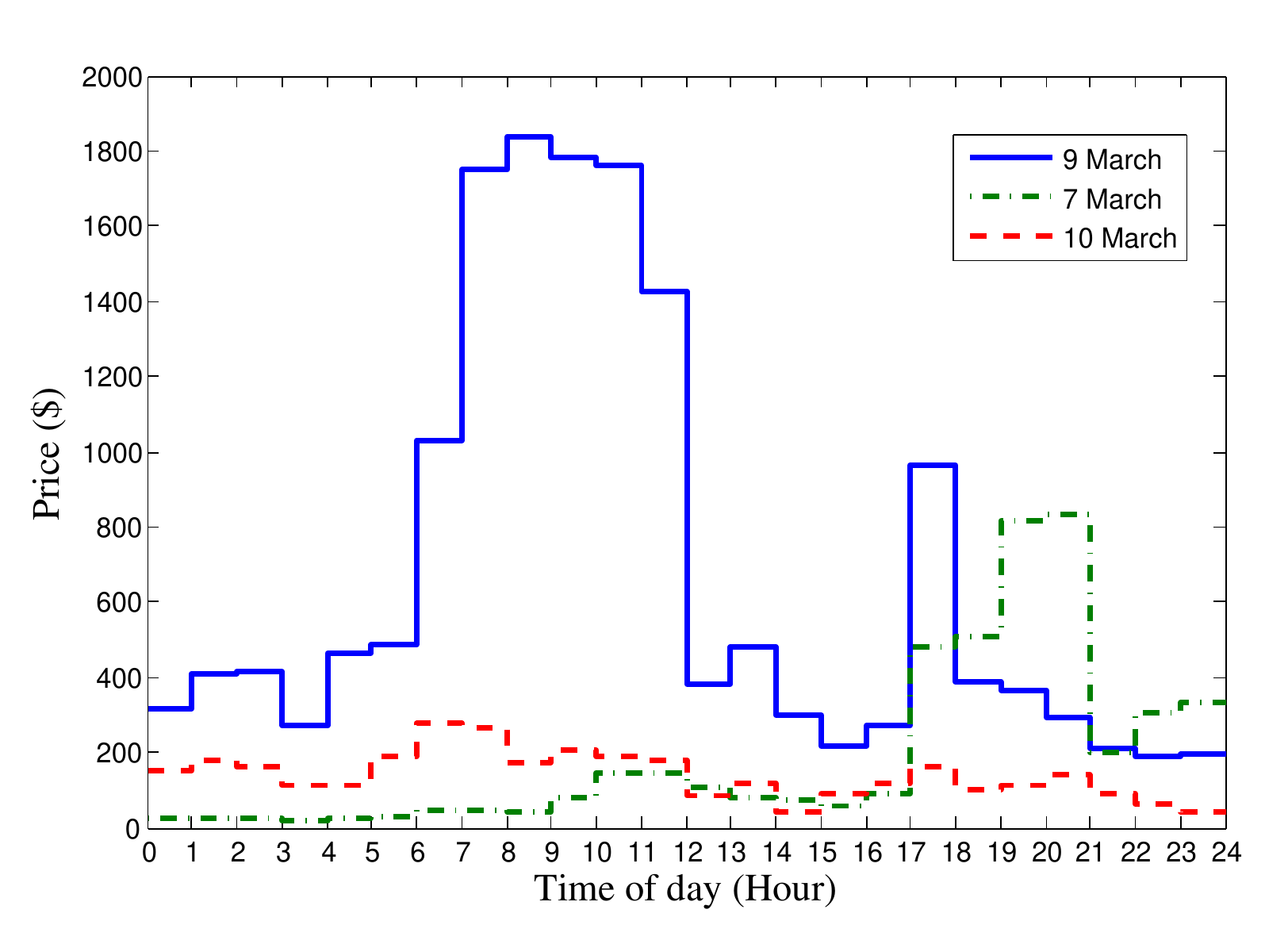}
	\caption{Volatility in the spot market for three consecutive days.} 
	\label{vol}
\end{figure} 

Fig. \ref{price} shows the DA and RT prices for 9 March 2014 as an example. We selected this day since it had the highest unexpected peak in RT prices throughout that year in PJM interconnection and can be assumed as what we earlier referred to as a black swan event in the market. As it can be observed, the RT price has a substantial peak at 9.00 A.M. around which the price is still unexpectedly high for 5 hours. The retailer proceeds for demand altering program upon receiving the RT pricing information according to \textbf{P2}. Fig. \ref{vol} shows that RT prices are totally different for one day before and one day after that day.  



\begin{figure} 
	\centering
	\includegraphics[width=\columnwidth,height=3.6in]{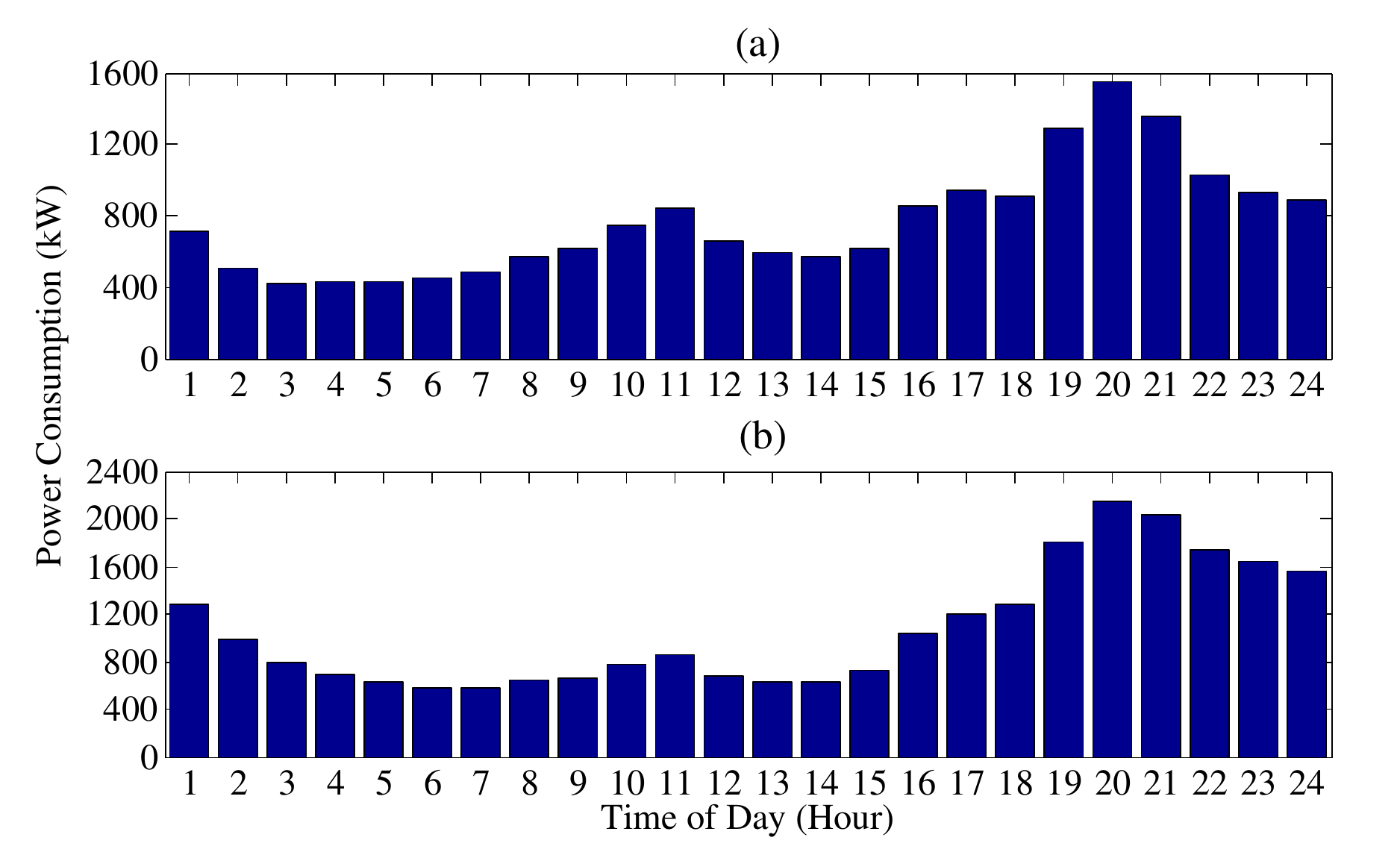}
	\caption[]{Electricity demand profile from (a) normal household appliances, i.e., without PEVs and (b) the overall electricity demand profile when users use PEVs with different usage patterns based on NHTS data.} 
	\label{infev}
\end{figure}

Next, we examine the DR scheme introduced in Algorithm 1. Fig. \ref{infev} shows the electricity demand profile from only typical household appliances, i.e., without PEVs and the overall electricity demand profile when users use PEVs with different usage patterns based on NHTS data.
 
To model the power purchased from the DA market on each day, we assume that the amount of power cleared at each hour of the following day in the DA market depends on the required power at that time slot and also includes some randomness as follows:

\begin{align}
d_{\text{cleared}}(t)=d_{\text{required}}(t)+u, \quad \forall t,
\label{dadem}
\end{align}

where $u$'s are independent and identically distributed (i.i.d.) chosen from a uniform distribution with the following representation:

\begin{align}
u \thicksim U[-0.2\times d_{\text{required}}(t),0.2\times d_{\text{required}}(t)) \quad \forall t.
\end{align}

We used a uniform distribution here as it has the highest entropy among bounded distributions and gives the highest uncertainty \cite{cover2012elements}.

\begin{figure} 
	\centering
	\includegraphics[width=\columnwidth,height=2.4in]{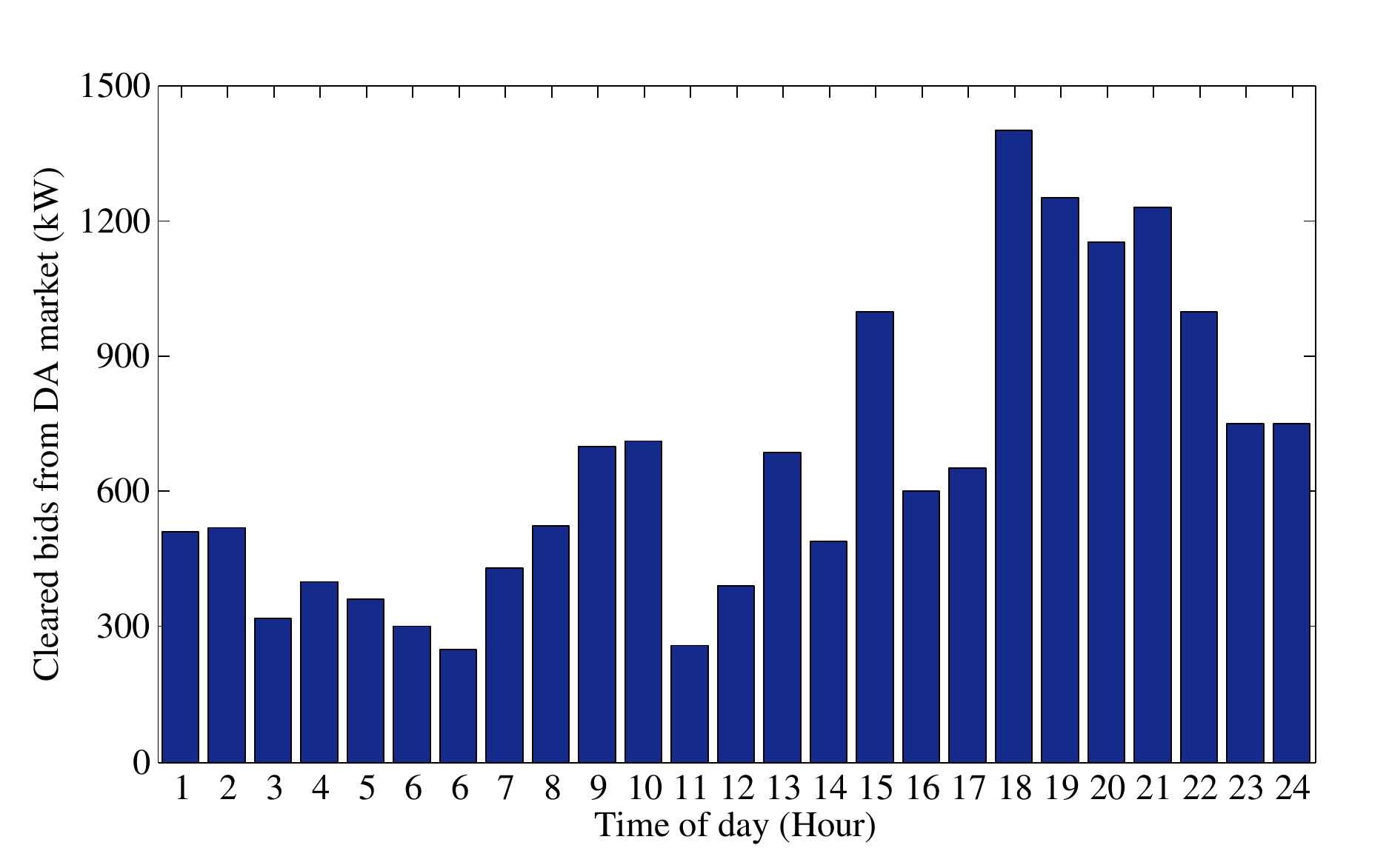}
	\caption[]{An assumed electricity profile purchased by the retailer from the DA market by the bids that could be cleared.} 
	\label{dap}\vspace{-2em}
\end{figure}



\begin{figure} 
	\centering
	\includegraphics[width=\columnwidth,height=3.6in]{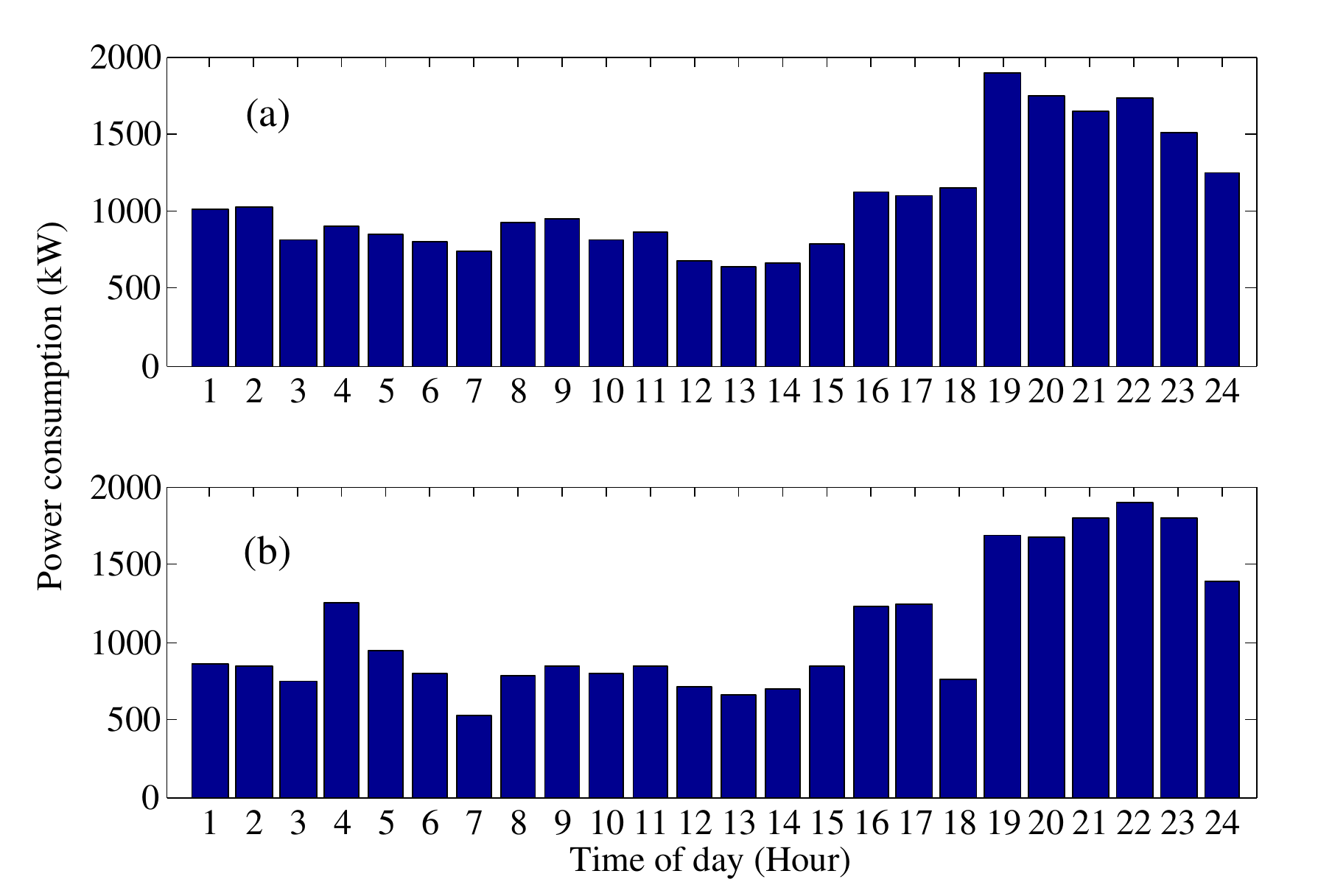}
	\caption[]{Load profiles after using Algorithm 1 for (a) just shaping the demand \textbf{P1} and (b) both shaping and altering the demand, Algorithm 1, on 9 March 2014.} 
	\label{shre}
\end{figure}

Fig. \ref{dap} shows the electricity profile purchased by the retailer from the DA market according to (\ref{dadem}) for different hours of a particular day. The results for shaping the DA demand profile only and joint shaping and altering the demand, Algorithm 1, is depicted in Fig. \ref{shre} \footnotemark        
\footnotetext[1]{The simulation was done using $\text{MATLAB}^{\circledR}$ 7.12.0. The mean simulation time was 172 seconds on a PC with intel i7-2670QM 2.2 GHz CPU, 6 GB RAM and $\text{Windows}^{\circledR}$7 OS.} . For this, we assumed $\lambda=0.5$ in \ref{P2}. In our simulations, convergence has been attained only after one single iteration of updating \textit{all} users' electricity demand profiles for both \textbf{P1} and \textbf{P2} in Algorithm 1.  


\begin{table}
	
	\renewcommand{\arraystretch}{8.5}
	\caption{\textlcsc{Overall energy procurement costs for the retailer on 9 March 2014}}
	\centering
	\resizebox{\columnwidth}{!}{  
		\begin{tabular}{|c|c|c|c|c|}
			\hline
			\Huge{Case} & \Huge{Overall ideal cost (\$)} & 
			\Huge{Overall real cost (\$)}
			&
			\Huge{Overall cost after \textbf{P1}} (\$)
			& \Huge{Overall cost after \textbf{P2} (\$)} \\ 
			\hline
			\Huge{1} & \Huge{\$674.4} &  \Huge{\$2,808.9} & \Huge{N/A} & \Huge{N/A}\\
			\hline
			\Huge{2}  &  \Huge{\$920.4} & \Huge{\$5,865.1} & \Huge{\$4,775.6} & \Huge{\$4,308.9}\\
			\hline
		\end{tabular}\vspace{-2.5 em}
	}
	\label{tab1}
\end{table}

In Table \ref{tab1}, we compare the overall energy procurement costs on 9 March 2014 for the retailer in two cases: \textit{case 1)} when there is no PEV in the system and \textit{case 2)} when all the users possess PEVs with their respective usage patterns. 

\begin{table}
	\renewcommand{\arraystretch}{8.5}
	\caption{\textlcsc{Overall energy procurement costs for the retailer in year 2014}}
	\centering
	\resizebox{\columnwidth}{!}{  
		\begin{tabular}{|c|c|c|c|c|}
			\hline
			\Huge{Case} & \Huge{Overall ideal cost (\$)} & 
			\Huge{Overall real cost (\$)}
			&
			\Huge{Overall cost after \textbf{P1}} (\$)
			& \Huge{Overall cost after \textbf{P2} (\$)} \\ 
			\hline
			\Huge{1} & \Huge{\$246,156} &  \Huge{\$921,664.3} & \Huge{N/A} & \Huge{N/A}\\
			\hline
			\Huge{2}  &  \Huge{\$324,778.3} & \Huge{\$1,840,333.8} & \Huge{\$1,761,889.1} & \Huge{\$1,688,926.1}\\
			\hline
		\end{tabular}\vspace{-2.5 em}
	}
	\label{tab2}
\end{table}

It can be noticed that in case 1, if the retailer could be absolutely successful in bidding to the DA market, i.e., there would not be any need to purchase electricity from the RT market (ideal occasion), the overall cost is only \$674.4, for the pricing shown in Fig. \ref{price}. In a more realistic case, when the retailer's bidding to the DA market is assumed to be according Fig. \ref{dap}, and the retailer is required to balance the load, the overall cost is remarkably higher. Obviously, demand shaping and demand altering in this case are not applicable (N/A) as there is no PEV and hence no power demand elasticity in the system.  

For the second case, when users possess PEVs, for the ideal bidding, the overall cost increases by almost 37\% to supply electricity to the PEVs whereas for the realistic bidding it becomes more than double. When demand shaping \textbf{P1} is employed, this overall cost reduces by around 18\%. Furthermore, when joint demand shaping and altering in Algorithm 1 is used, cost decreases further by almost 10\%.  

In this case, we also assumed that the retailer is allowed to sell back its extra load purchased earlier from the DA market to the RT market at the same RT prices in the RT market. 

In Table \ref{tab2}, we evaluate the proposed technique in Algorithm 1 over the whole year of 2014 according to the PJM interconnection data. In this case, as stated in the algorithm, the retailer proceeds for minimizing its RT demand if $p^r_t>\gamma_K(t)$ and its maximize its purchase from RT market if $p^r_t<\gamma_K(t)$. It is observed that \$151,407.7 can be saved in the whole year for the energy procurement cost of the system by using the proposed algorithm.


\section{Conclusion and Future Work} \label{Con}

In this paper, we proposed a fast converging and decentralized algorithm for managing V2G enabled PEVs' electricity assignments (charging and discharging) to lower the overall electricity procurement cost for an electricity retailer. Our proposed algorithm uses demand shaping and demand altering for the DA and the RT markets. In particular, when the power system has high penetration of intermittent energy resources, demand altering is crucial due to likely contingencies and hence more RT price fluctuations. In our simulations' results, we considered the pricing data in PJM interconnection electricity market for the year 2014. We showed that significant overall cost savings (up to \$151,407.7) for a retailer bidding to this electricity market could be achieved by using our proposed algorithm throughout the year. This allows the retailer to offer better deals to the customers and expand its market capacity and customers can enjoy lower electricity bills as well.  

The work presented in this paper can be extended in various ways. First, recently, PEVs' battery storage capacity has been expanded up to 60-85 kWh, e.g., for Tesla model S and model X \cite{Tesla}, which can provide much more elasticity for demand shaping and demand altering and hence reduce the electricity costs further. Second, our emphasis in this paper was on active power, this work can be extended to reactive power as well.         

%
%

\bibliographystyle{IEEEtran} 
\bibliography{IEEEabrv,myBIB}

\begin{footnotesize}
	\begin{IEEEbiography}
		[{\includegraphics[width=1in,height=1.25in,clip,keepaspectratio]{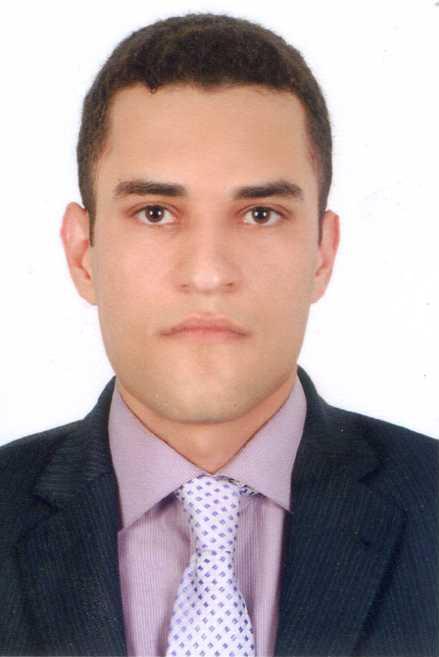}}]%
		{Farshad Rassaei}
		(S'10) has received his B.Sc. and M.Sc. degrees with first class honours both from the department of electrical and computer engineering, Shiraz University, Shiraz, Iran, in 2009 and 2012, respectively. He is currently working toward the Ph.D. degree at the department of electrical and computer engineering, National University of Singapore (NUS), Singapore. He is a recipient of Singapore international graduate award (SINGA). He has also served as the President of the fifth NUS ECE Graduate Student Symposium (GSS2015). His research interests include communication systems, signal processing, optimization theory, game theory and smart grids. 
	\end{IEEEbiography}
	\begin{IEEEbiography}[{\includegraphics[width=1in,height=1.25in,clip,keepaspectratio]{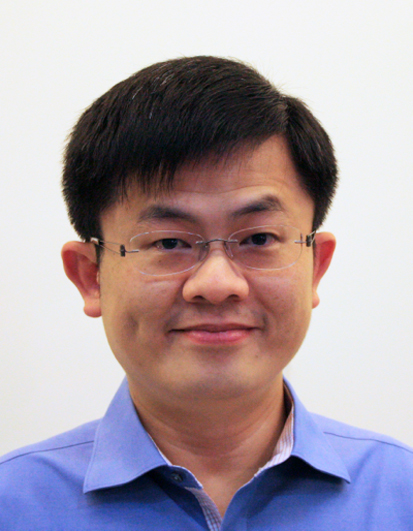}}]{Wee-Seng Soh} (S'95-M'04) received the B.Eng. (Hons) and M.Eng. degrees in electrical engineering from the National University of Singapore (NUS) in 1996 and 1998, respectively. In 1998, he was awarded the Overseas Graduate Scholarship by the National University of Singapore to study at Carnegie Mellon University, Pittsburgh, PA, where he received the Ph.D. degree in electrical and computer engineering in 2003. Since 2004, he has been with the Department of Electrical and Computer Engineering, National University of Singapore, where he is currently an Associate Professor. Prior to joining NUS, he was a Postdoctoral Research Fellow in the Electrical Engineering and Computer Science Department, University of Michigan, Ann Arbor. He has served on the Technical Program Committees (TPC) of over 20 conferences, and has also served as a TPC co-chair for the 2008 IEEE International Conference on Communication Systems (ICCS), the 2011 IEEE International Workshop on Underwater Networks (WUnderNet), and the 2013 IEEE International Conference on Networks (ICON). He is currently serving as an Area Editor of Computer Communications (Elsevier), and also, Pervasive and Mobile Computing (Elsevier). His current research interests are in wireless networks, underwater networks, indoor tracking/localization techniques, and satellite systems.
	\end{IEEEbiography}
	\begin{IEEEbiography}[{\includegraphics[width=1in,height=1.25in,clip,keepaspectratio]{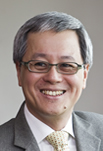}}]{Kee Chaing (KC) Chua} received a PhD degree in Electrical Engineering from the University of Auckland, New Zealand in 1990 and joined the Department of Electrical Engineering at the National University of Singapore (NUS) as a Lecturer.  He is now a Professor in the Department of Electrical \& Computer Engineering at NUS.  He served as the Faculty of Engineering's Vice Dean for Research twice, from 2003 to 2006 and from 2008 to 2009.  From 1995 to 2000, he was seconded to the Centre for Wireless Communications (now the Institute for Infocomm Research), a national telecommunication R\&D centre funded by the Singapore Agency for Science, Technology and Research as its Deputy Director.  From 2001 to 2003, he was on leave of absence from NUS to work at Siemens Singapore where he was the Founding Head of the Mobile Core R\&D Department funded by Siemens' ICM Group.  From 2006 to 2008, he was seconded to the National Research Foundation as a Director.  He was appointed Dean of the Faculty of Engineering at NUS in July 2014, after serving as Head of its Department of Electrical \& Computer Engineering since November 2009.  He chaired the World Economic Forum's Global Agenda Council on Robotics and Smart Devices in 2011 and spoke on the role of robotics and smart devices in shaping new models of development at the World Economic Forum in Davos in January 2012.  He is a recipient of an IEEE 3rd Millennium Medal, a Fellow of the Singapore Academy of Engineering, and a Fellow of the Institution of Engineers, Singapore.
	\end{IEEEbiography}
\end{footnotesize}

\end{document}